**A Dynamical Systems Framework for Resilience in Ecology**

Katherine J. Meyer

University of Minnesota

127 Vincent Hall
206 Church St. SE
Minneapolis, MN 55455

**Abstract**

Rising interest in the resilience of ecological systems has spawned diverse interpretations of the term's precise meaning. This paper classifies and explores definitions of resilience from the ecological literature using a dynamical systems framework. A model consisting of ordinary differential equations is assumed to represent the ecological system. The question "resilience of what to what?" posed by Carpenter et al. [2001] informs two broad categories of definitions, based on resilience to state variable perturbations and to parameter changes, respectively. Definitions of resilience to state variable perturbations include measures of basin size (relevant to one-time perturbations) and of basin steepness (relevant to repeated perturbations). Resilience to parameter changes can be quantified by viewing parameters as state variables, but could also take into account reversibility. The dynamical systems viewpoint yields the following key insights: (i) system properties that confer resistance to state variable perturbations can differ from those that promote recovery from them, (ii) recovery rates, oft deemed the purview of "engineering resilience," do matter for ecological resilience, and (iii) resilience to repeated state variable perturbations correlates with resilience to parameter changes due to the critical slowing down phenomenon.

**Keywords:** resilience, dynamical systems, ecology

# 1. Introduction

Ecological systems as diverse as coral reefs, grasslands and tropical forests can gravitate toward alternative stable states (Folke et al. [2004]), which provide differing flows of ecosystem services to humans. Our simultaneous dependence on and disruptions to ecosystems make a strong case for careful natural resource management. Promoting ecosystem *resilience*—broadly defined as the ability to absorb disturbance while maintaining the same basic structure and function—has emerged as an important alternative to the maximum sustainable yield paradigm (Walker and Salt [2006]). As the use of "resilience" in the sustainability vernacular has grown, so has the variety of interpretations of its meaning. To be most useful, definitions of resilience should be flexible enough for wide application but clear enough to avoid confusion.

Other authors have proposed frameworks to clarify the meaning of resilience. Holling [1996] distinguished between "ecological resilience," concerned with staying near one of multiple alternative stable states, and "engineering resilience," concerned with the return time to equilibrium following a perturbation. By asking "resilience *of what to what*?" Carpenter and colleagues [2001] added an important lens to the discussion of resilience.

This paper builds on the work described above, further organizing ecological definitions of resilience using the dynamical systems framework of mathematics. Section 2 provides a summary of dynamical systems modeling as it applies to ecology. Next, sections 3 and 4 address the resilience *of* a system staying in the same basin of attraction, *to* state variable perturbations (section 3) and parameter changes (section 4). Finally, section 5 presents some considerations related to timescale, and Table 1 summarizes key aspects of resilience treated in the text.

The dynamical systems language of mathematics provides a useful framework for reconciling and clarifying resilience definitions. In applying resilience theory to real world systems, however, quantification is not always possible and qualitative understanding certainly yields useful results (Walker et al. [2004]). Quantitative analysis can, however, reveal effects that appear counterintuitive from a qualitative viewpoint. These include (i) the difference between the strength of recovery of a perturbed variable and the difficulty of causing the perturbation (section 3.1), (ii) the importance of recovery rates for even Holling's ecological resilience (section 3.3), and (iii) the interconnection between resilience to state variable and parameter changes (section 4). Rather than supplanting qualitative uses of resilience, I hope that this mathematical treatment will add clarity and insight to how we might promote ecosystem resilience.

## 2. Modeling outlook

### 2.1 State variables and parameters

For an excellent overview of dynamical systems modeling in ecology, see Beisner et al. [2003b]. Any mathematical model must make assumptions about what quantities matter for understanding the question at hand. Dynamical systems models include two types of quantities: state variables and parameters. State variables are best suited for quantities that change relatively rapidly through feedbacks with one another, such as population densities. In a continuous dynamical system—the type considered here—differential equations describe the rates of change in each state variable as functions of the state variables and parameters. Parameters thus influence the state variable evolution but don't change dynamically themselves; they are often appropriate for slow-changing quantities such as environmental factors. For example, a fisheries model might include trout population as a state variable and environmental carrying capacity as a parameter. The modeler designates state variables and parameters to suit the questions at hand.

## 2.2 The stability landscape

The system behavior that arises from a model's differential equations can be explored mathematically via analytic, geometric, or numerical approaches. Ecologists often use the stability landscape as a geometric technique for visualizing dynamics (Figure 1). Horizontal positions of an imagined ball on the stability landscape encode the value of state variable(s); feedbacks and parameters shape the topography of the surface. In mathematical parlance, the surface of a stability landscape represents a potential function for the vector field that describes the rates of change in state variables. The system evolves on the stability landscape under "negative gradient flow" by rolling in the steepest downhill direction at a rate proportional to the slope.

Flat spots on the stability landscape represent equilibria for the system. These may be stable (**A**, **C**) or unstable (**B**), depending on the local slope of the landscape. Alternative stable states (**A** and **C**), which arise from nonlinear feedbacks between state variables, appear as distinct troughs or wells on the stability landscape. For each attracting stable state, the collection of positions that flow to the stable state under the feedbacks of the system is called its *basin of attraction* or *domain of attraction*.

Stability landscapes have important limitations. First, the metaphor of a ball rolling on the surface can be interpreted too literally. Rather than rolling according to the force of gravity and passing through flat regions using momentum, the ball travels with velocity dictated at each moment by the slope of the landscape in the most downhill direction. As the slope of the stability landscape decreases to zero near a stable equilibrium, the ball slows and approaches the equilibrium asymptotically. Second, systems of differential equations with two or more state variables may lack a corresponding surface for negative gradient flow due to rotational tendencies in the vector field. For example, a system that cycles periodically in a circular pattern can't perpetually flow downhill on a surface. Although a stability landscape still exists for these systems, it describes the system's behavior more loosely: the ball can't roll uphill, but it can go downhill (not necessarily parallel to the gradient) or move along a level curve. Lastly, stability landscapes are difficult to visualize for systems with three or more state variables. Given both their utility and limitations, let us consider stability landscapes a helpful cognitive device for discussing resilience, interpreting them carefully.

### 2.3 Ignored quantities

In order to remain tractable for analysis, a model cannot include everything that might influence the question at hand. Good models give insight into the mechanisms that drive a system, in part because of simplifying assumptions that highlight key processes. In the context of resilience, unexpected disturbances arise from processes that aren't represented in the model. We turn to state variable disturbances in section 3, followed by parameter changes in section 4.

## 3. State variable perturbations

### 3.1 Source of state variable perturbations

Before discussing resilience to state variable perturbations, we consider their origins. State variable perturbations arise from variations in quantities and processes not included in a model (Beisner et al. [2003b]), called "disturbances" (Pimm [1984]) or "external forces" (Carpenter et al. [2001]). The term *resistance* describes the relationship between external forces and the resulting state variable perturbation (Levin and Lubchenco [2008], Carpenter et al. [2001]). A highly resistant system will register only small changes to state variables under considerable external forcing (Table 1, row 1). For example, if a vessel of liquid is our system and the liquid's temperature the state variable, a liquid with high heat capacity has high resistance to external heating or cooling disturbances.  Unlike resilience, which asks whether a system can dynamically recover from a state variable perturbation, resistance addresses whether significant perturbations will occur in the first place. Together, resistance and resilience to perturbations determine the persistence of a system state through time (Carpenter et al. [2001]).

The description of resistance given by Walker et al. [2004] as the slope of a basin of attraction warrants scrutiny. They claim that "steep sides imply greater perturbations or management efforts are needed to change the state of the system, i.e., its position within the basin" (p. 3). Assuming that the shape of the basin of attraction represents the potential function for the system's vector field, this statement could be misleading. Because changes outside of the modeling framework cause state variable perturbations, they are not necessarily influenced by the feedbacks that control the state variables, represented by the slope. This is true especially when timescales of disturbance and recovery differ.  For example, in the case of a fire disturbance to vegetation, the properties of vegetation that minimize acute fire damage are distinct from the properties that promote long-term regrowth following the disturbance. However, by affecting the post-perturbation response, basin slope could certainly affect resilience to repeated state variable perturbations (section 3.3).

### 3.2 Resilience to a one-time state variable perturbation

Standard mathematical definitions of stability such as Lyapunov stability fail to capture the notion of resilience to state variable perturbations. An attractor in a dynamical system is called *Lyapunov stable*  if for any neighborhood $U$ of the attractor, there's a neighborhood $V$ such that the system will surely stay within *U* as time goes forward, so long as it starts within *V*. Importantly, this notion of stability doesn't specify the size of the neighborhood *V* that guarantees staying in *U*. Because of this, Lyapunov stability does not indicate whether an attractor will be resilient to state variable perturbations of a certain magnitude and direction. Further characterizations of the domain of attraction are needed.

To quantify a domain of attraction's resilience to a single perturbation, various measures of basin size may be relevant. Roughly speaking, bigger basins require larger perturbations for escape. Possible quantities include the width (1D) or area (2D) of the basin (Figure 2a,b; Table 1, row 2). Together, these types of measurements have been termed "flexibility" (Levin and Lubchenco [2008]) or "latitude" (Walker et al. [2004]).

However, a basin could extend infinitely in one direction yet finitely in others, yielding infinite width or area but limited resilience. In such cases, the distance from a desired stable state to a threshold between basins (Beisner et al. [2003a]) offers more insight into the size of perturbations from which the system could recover (Figure 2c; Table 1, row 3). This simple example reveals that even in a one-dimensional state space, direction matters. One could specify the distance to a threshold for a particular perturbation direction of interest, or take the minimum over set of anticipated directions (Figure 2d).

The measures described thus far characterize an entire basin of attraction, independent of the system's current position within it. They are probably most informative when the system tracks the equilibrium closely. On the other hand, repeated perturbations could cause a system to wander considerably about the domain of attraction, making the distance from the current state to a threshold (termed "precariousness" by Walker et al. [2004]) a relevant measure of resilience to further perturbations (Table 1, row 4). We take up the topic of resilience to repeated state-space perturbations next.

### 3.3 Resilience to repeated (smaller) state-space perturbations

When a domain of attraction's resilience to not just one but repeated state variable perturbations matters, the strength and direction of the vector field within the domain becomes important. Large-magnitude vector fields directed away from thresholds will help the system recover back towards the attracting equilibrium between perturbations, making escape more difficult. Various names have been given to describe this property of a vector field, including "intensity of attraction" (McGehee [1988]), "slope" (Beisner et al. [2003b]), "recovery" (Levin and Lubchenco [2008]), and "resistance" (Walker et al. [2004]).

Fully characterizing the magnitudes and frequencies of state variable perturbations that a basin of attraction can absorb presents an interesting mathematical task. A key factor controlling resilience under a given perturbation frequency is how far the system can recover back towards equilibrium between perturbations. This maximum recovery, which could occur anywhere in the basin, depends on the strength of the vector field over that region (Table 1, row 5).

Restricting attention to steepness *near a stable equilibrium* yields Pimm's [1984] definition of resilience, namely "how fast the variables return towards their equilibrium following a perturbation." A linear system $x'=Ax$ recovers to *1/e* of a perturbation's magnitude in roughly *-1/Re(λ)* time units, where λ is the weakest eigenvalue of *A* (Neubert and Caswell, [1997]). Pimm defines resilience as the reciprocal of this characteristic return time. Holling [1996] calls this definition "engineering resilience," and argues against its use in ecology because it takes a linearized view of systems, neglecting the possibility of alternative stable states. Nonetheless, return rates do matter for Holling's ecological resilience in the context of repeated state variable perturbations . Simply put, "steep parts" help the system recover towards equilibrium, but could occur anywhere in the domain of attraction.

In contrast to other definitions, Ives and Carpenter [2007] quantified Holling's ecological resilience as the rate of repulsion at an unstable equilibrium between stable states. Repulsion rates near an unstable equilibrium provide an incomplete picture of resilience of a domain of attraction for the same reason that return rates near a stable equilibrium do.  As is shown for a simple 1D example in Figure 3, the region of the basin of attraction that most thwarts escape under repeated small perturbations may occur somewhere between the stable and unstable equilibria. Imagine that the system starts at rest at the stable equilibrium on the left, but the state variable is perturbed to the right repeatedly by external forces. Near the stable equilibrium, the state variable will flow to the left only weakly during the recovery time between perturbations. The state may continue to move to the right in the presence of the disturbance.  Once it reaches the steep section, however, the leftwards flow between disturbances will increase so much that it may fully balance the disturbances, preventing the state variable from moving any further towards the basin threshold.

## 4. Resilience and Parameter Changes

### 4.1 Parameters affect state variable resilience

Events outside the modeling framework can alter the values of not only state variables but also parameters. Importantly, parameter changes alter the resilience of the system to state variable perturbations, as defined above (Beisner et al. [2003b]). Figure 4a illustrates for the simple system $x' = \mu + x - x^3$ how changes in the parameter $\mu$ gradually reduce the size and steepness of the basin of attraction on the left for the state variable $x$.

### 4.2 Resilience to parameter changes

In addition to measuring how parameters affect resilience to state variable perturbations, we may also consider a system's resilience *to parameter changes themselves*. Because the distinction between state variables and parameters is arbitrary, one approach simply reclassifies parameters as state variables with zero rate of change and then applies the definitions of resilience to state variable perturbations (Beisner et al. [2003b]). Figure 4b shows the two-dimensional state space that results from making $\mu$ a state variable in the system $x' = \mu + x - x^3$. Since dynamics evolve in only the vertical direction, the domains of attraction for the lower and upper branches of equilibria are the dark blue and light orange regions, respectively. The horizontal distance (i) from the current $\mu$ value to the bifurcation value $\mu$* yields Ives and Carpenter's [2007] definition of resilience to parameter changes, namely the "magnitude of the press perturbation before the stable equilibrium bifurcates" (Table 1, row 6). Precariousness (ii)—the distance from a point in $(x, \mu)$ space to the other basin of attraction— could also be considered (Table 1, row 7).

Note that even when treated like state variables, parameters don't dynamically change. With parameter recovery rates of zero, defining resilience to repeated parameter perturbations doesn't make sense. A system's resilience to repeated *state variable* perturbations, however, may reveal something about its resilience to *parameter* changes. The phenomenon of reduced recovery rates near a bifurcation—termed critical slowing down—has been proposed as a way to predict a catastrophic regime shift (Scheffer et al. [2009], Wissel [1984]). Because basins tend to get shallower as they approach a bifurcation, resilience to repeated state variable perturbations correlates with the magnitude of parameter changes that can be sustained without losing the basin of attraction.

When parameters *do* change enough to move the system to a fundamentally different domain of attraction, questions of resilience turn to the reversibility of the change. Following the lead of Carpenter et al. [1999], we can classify state shifts caused by parameter changes into three types: reversible, hysteretic, and irreversible. Figure 5 illustrates the key properties of each type of parameter change.

A regime shift caused by a *reversible* parameter change can be undone by simply returning the parameter to its initial value (Figure 5a). A *hysteretic* system can also recover its initial state following a regime shift, but restoring the initial parameter value does not suffice because the system gets stuck at the undesired alternative stable state; further parameter change is needed to precipitate recovery (Figure 5b). An *irreversible* parameter change arises when physical constraints prohibit parameter values that would restore the desired state in a hysteretic system (Figure 5c). For example, in the lake model of Carpenter et al. [1999], phosphorous inputs from the watershed (a parameter) control the water quality state variable. Irreversible lakes cannot return to the clear water state via manipulation of phosphorous loads alone because this would require negative phosphorous loading. Importantly, "irreversibility" depends on the manipulations one allows; dredging lake bottoms to remove sediment phosphorus can alter the shape of the hysteretic curve and restore clear water to a lake deemed irreversible with respect to phosphorous loading. In sum, the concept of reversibility broadens potential definitions of parameter resilience to include not only the difficulty of reaching bifurcations but also the possibility of recovering from them (Table 1, row 8).

## 5. Timescale

The qualitative definition of resilience—a system's ability to absorb disturbance while maintaining the same basic structure and function—leaves timescale open for interpretation. Holling took a long-term view when he articulated the concept of resilience in [1973]. Among resilient systems, he included ecosystems such as forests that cycle through *different* basins of attraction over the course of decades. We could even imagine systems that cycle on a geologic timescale through glacial periods. In this perspective, variation is part of the structure and function of the ecosystem, and it may take long time to return to a particular basin after a disturbance. Since Holling's seminal work, he has joined other resilience thinkers in defining resilience as the ability to stay in the *same* basin of attraction. Nonetheless, many resilience thinkers continue to take a long-term view: slow returns within a basin, such as the decadal recovery of a coral reef following a hurricane, still reflect system resilience (Walker and Salt [2006]).

When humans rely on a particular system for ecosystem services, natural resource managers may desire more than long-term resilience. Knowing that an ecosystem service will be restored within the next glacial cycle does little for today's populations. Rouge et al. [2013] emphasized the importance of defining a decision-making time horizon when they explored resilience definitions from a complex systems engineering perspective. Under this view, a system can be resilient under a 50 year time horizon but not a 20 year time horizon. These interests in getting back to a desired state either eventually or on a management-decided timeframe reflect a tension between letting nature take its course and interfering for shorter-term human benefit.

Timescale also plays an important role in the theory of dynamical systems. Focusing on the long-term behavior of dynamical systems spurred great theoretical advances during the 20[th] century. However, short-term system behavior can differ significantly from long term behavior (Neubert and Caswell [1996]), and interactions between short term system feedbacks and external perturbations can have important consequences for resilience. Developing dynamical systems theory to understand these short term interactions presents an exciting area for future research. (Zeeman [2014]).

## 6. Summary and Conclusions

Table 1 summarizes the major indicators of resilience described in this paper, and proposes potential mathematical definitions of each. Linking resilience to a dynamical systems framework helps organize the varied concepts of resilience and clarifies their connections. For example, this framework revealed that two uses of the term "resistance" that actually refer to different concepts: the difficulty of perturbing a state variable by external force (Levin and Lubchenco [2008], Carpenter et al. [2001]) and the strength of feedbacks that counteract perturbations once they occur (Walker et al. [2004]). Dynamical systems also show that return rates— despite their unpopularity in Holling's ecological-resilience camp—play an important role in determining resilience to repeated perturbations when they are considered over the entire basin of attraction and not simply at equilibria. Rather than aiming for one definition of resilience, we might embrace the variety of indicators available. When assessing the resilience of a real-world natural system, context could determine which indicators matter most.

Although many of the indicators in Table 1 translate neatly into a mathematical formulation, resilience to repeated state variable perturbations and reversibility of parameter changes don't. These two aspects of resilience appear ripe for contributions from mathematics research. Even for real-world systems that can't be quantified to the precision suggested by dynamical equations, research models could provide insight into possible phenomena and key variables to measure.  Investigating these aspects of resilience could thus advance both mathematical theory and sustainability science.

## Acknowledgements

Thank you to Allison Shaw, Kelsey Vitense, Lauren Sullivan, Ranjan Muthukrishnan, and Forest Isbell. Our conversations have shaped my understanding of resilience in ecology.

**Table 1** Resilience Indicators

| Resilience *to what*? | | Indicators | Mathematical definitions |
|---|---|---|---|
| **State variable perturbations (section 3)** | single | resistance | $\dfrac{\text{external disturbance}}{\text{resulting perturbation}}$ |
| | | basin area (2D) or width (1D) | $n$-dimensional volume of domain of attraction for an attractor in $n$-dimensional state space |
| | | distance from stable equilibrium to basin threshold | $\inf_s |\mathbf{a} - s|$ where $\mathbf{a}$ is stable equilibrium, $s$ in separatrix |
| | | distance from current position to basin threshold | $\inf_s |\mathbf{x} - s|$ where $\mathbf{x}$ is current state, $s$ in separatrix |
| | repeated | basin "slope" | *intensity of attraction* (McGehee 1988) *research area* |
| **Parameter changes (section 4)** | | distance from current parameter value(s) to bifurcation values | $\inf_{\mu^*} |\mu - \mu^*|$ where $\mu$ is current parameter value, $\mu^* \in \{\text{bifurcation values}\}$ |
| | | distance to threshold when parameter treated as state variable | $\inf_s |\mathbf{x} - s|$ where $\mathbf{x}$ is current state, $s \in$ separatrix in expanded state space |
| | | reversible / hysteretic / irreversible | *research area* |

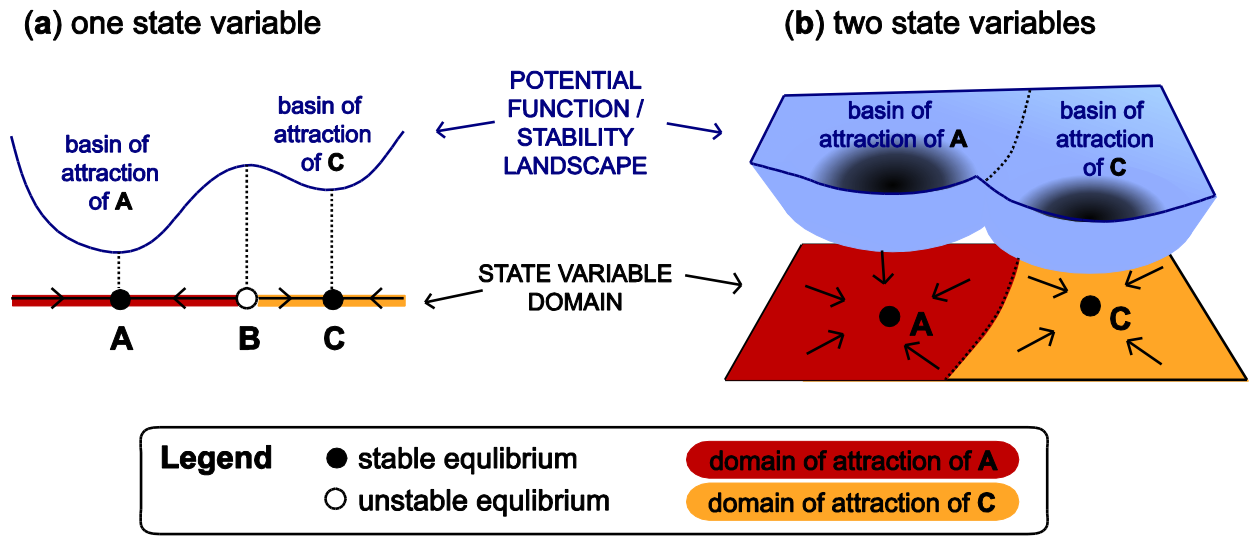

**Figure 1** Key terms related to stability landscapes.

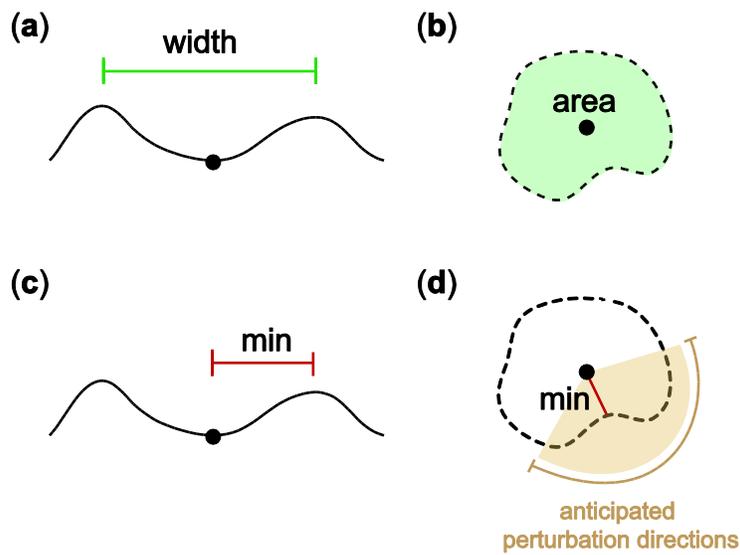

**Figure 2** Measurements of resilience to a single state variable perturbation include (**a**) width or (**b**) area of a basin of attraction, and (**c, d**) the minimum distance from an equilibrium to a threshold.

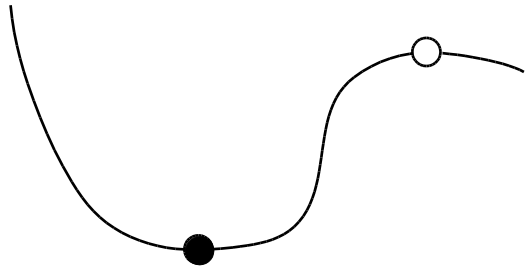

**Figure 3** Steep basin sections may occur far from equilibria. Stable (black) and unstable (white) equilibria are shown on the potential function for a one dimensional system.

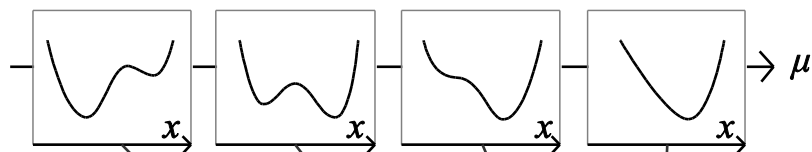

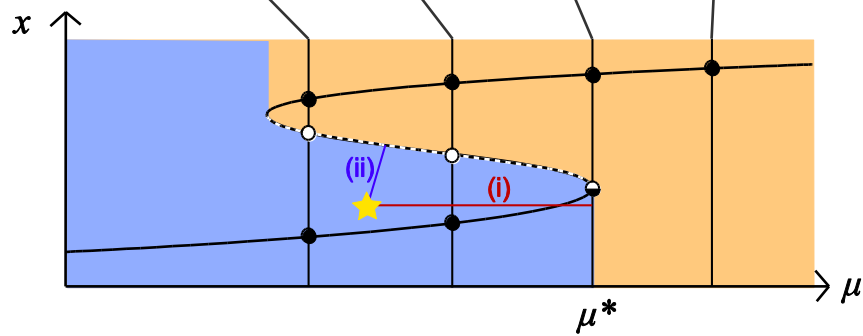

**Figure 4** Resilience and parameter changes. (**a**) As the parameter $\mu$ increases, the resilience of the left basin—measured either as distance from stable equilibrium to threshold or as steepest slope—decreases. (**b**) Each example in (a) corresponds to a vertical slice in $\mu$-$x$ state space. The endogenous evolution of the system occurs in the vertical direction only, and with blue and orange basins of attraction for the lower and upper branches of equilibria, respectively.

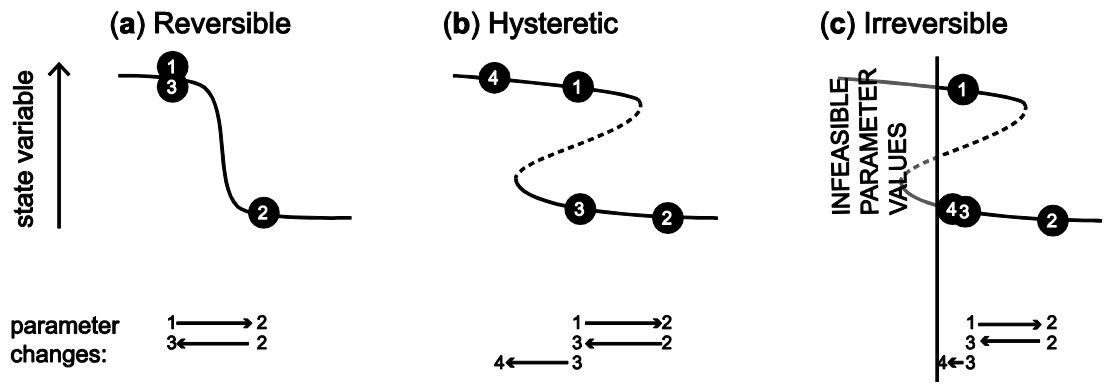

**Figure 5** Reversible, hysteretic, and irreversible parameter changes. In these bifurcation diagrams, solid and dashed lines indicate stable and unstable equilibria, respectively, in the space of parameter (horizontal axis) and state variable (vertical axis). Numbered circles indicate the sequential position of the system under the parameter changes indicated below. In each case, an initial parameter change shifts the system from a high to low equilibrium state. System (**a**) is reversible because reversing the parameter change also undoes the associated state change. System (**b**) is hysteretic because further parameter change beyond a simple reversal is necessary to restore the state variable to the upper branch of equilibria. System (**c**) is irreversible because the parameter changes necessary to restore the state variable are physically infeasible.